\documentclass[12pt]{article}
\usepackage[dvips]{graphicx}
\usepackage{amsmath,amssymb}
\usepackage{bm}

\title{Brief note on Riemann hypothesis II}
\author{}

\begin{document}
{\pagestyle{empty}
\rightline{November 2011}
\rightline{~~~~~~~~}
\vskip 1cm
\centerline{\large \bf A Brief Note on the Riemann hypothesis II}
\vskip 1cm
\centerline
{{Minoru Fujimoto
 } and
 {Kunihiko Uehara
 }
}
\vskip 1cm
\centerline{\it ${}^1$Seika Science Research Laboratory,
Seika-cho, Kyoto 619-0237, Japan}
\centerline{\it ${}^2$Department of Physics, Tezukayama University,
Nara 631-8501, Japan}
\vskip 2cm

\centerline{\bf Abstract}
\vskip 0.2in
  We have dealt with the Euler's alternating series 
of the Riemann zeta function 
to define a regularized ratio appeared in the functional equation
even in the critical strip and 
showed some evidence to indicate the hypothesis. 
We briefly review the essential points and 
we also define a finite ratio in the functional equation 
from divergent quantities in this note. 

\vskip 0.4cm\noindent
MSC number(s): 11M06, 11M26

\noindent
PACS number(s): 02.30.-f, 02.30.Gp, 05.40.-a

\hfil
}
\setcounter{equation}{0}
\addtocounter{section}{0}
\hspace{\parindent}

  Until now regularizations by way of the zeta function have been successful in 
some physical problems\cite{Kunimasa,Leeuwen}, but the Riemann hypothesis itself has been remained to be proved. 
Recently we have shown the evidence using the Euler's alternating summation, 
which is finite even in the critical strip and seems to be essential to clarify the Riemann hypothesis, 
in the story of the finite quantities.\cite{Fujimoto5}

  Here we will briefly review the essential points and give other evidences, one of which is also related with 
the finite ratio appeared in the functional equation and 
another is in the story of the divergent quantities for the previous derivation. 


  The definition of the Riemann zeta function is 
\begin{equation}
  \zeta(z)=\lim_{n\to\infty}\zeta_n(z),\ 
  \zeta_n(z)\equiv\sum_{k=1}^n\frac{1}{k^z}
\label{e101r}
\end{equation}
for $\Re z>1$ \cite{Conrey}.
In this note we adopt a hat notation such as $\hat{\zeta}(z)$, 
which is well defined even in the critical strip $0<\Re z<1$, 
such as the Euler's alternating series
\begin{equation}
  \hat\zeta(z)=\frac{1}{1-2^{1-z}}\lim_{n\to\infty}\xi_n(z),\
  \xi_n(z)\equiv\sum_{k=1}^n\frac{(-1)^{k-1}}{k^z}.
\label{e102r}
\end{equation}
We often mention a hat notation as a ``regularized" form because 
a hat expression is defined by a subtraction 
an infinite number from a divergent quantity.

  In this note, we deal with the Euler's alternating series 
of the Riemann zeta function as (\ref{e102r}) 
to well-define even in the critical strip $0<\Re z<1$ 
and utilize the functional equation to indicate the hypothesis. 
Hereafter we are only interested in the region $\Re z\ge\frac{1}{2}$ 
for the Riemann zeta function, because the functional equation ensures 
the regularized nature of the zeta function for the other half plane 
$\Re z<\frac{1}{2}$.

\vskip 5mm
  There is a relation called the functional equation for the Riemann zeta function 
\begin{equation}
  \hat{\zeta}(z)=\hat{H}(z)\hat{\zeta}(1-z),
\label{e121r}
\end{equation}
where $\hat{H}(z)$ is given by $\displaystyle 2\Gamma(1-z)(2\pi)^{z-1}\sin\frac{\pi z}{2}$ 
and is not equal to zero for $\displaystyle{\frac{1}{2}\leq\Re{z}<1}$.
Hereafter we deal with $\hat{H}(z)$ as the infinite limit of $\hat{H}_n(z)$ defined by 
\begin{equation}
  \hat{H}_n(z)\equiv\frac{\hat{\zeta}_n(z)}{\hat{\zeta}_n(1-z)},
\label{e122r}
\end{equation}
where $\hat{\zeta}_n(z)$ is defined by
\begin{equation}
  \hat{\zeta}_n(z)\equiv\zeta_n(z)-\frac{n^{1-z}}{1-z}
\label{e123r}
\end{equation}
but as we will give a notice, we have to take care of substituting a zero for $z$ 
in the limit of $n\to\infty$. 

The relation between $\zeta_n(z)$ and $\xi_n(z)$ is special 
because the relation form itself conserves before and after the regularization 
as
\begin{eqnarray}
  \xi_{2n}(z)&=&\zeta_{2n}(z)-2^{1-z}\zeta_n(z),
\label{e202ar}\\
  \xi_{2n}(z)&=&\hat{\zeta}_{2n}(z)-2^{1-z}\hat{\zeta}_n(z),
\label{e202br}
\end{eqnarray}
where we used the relation (\ref{e123r}) 
and we do not use a hat notation for $\xi_{2n}(z)$ 
because it is already well-defined in the critical strip.
Adding the term $\hat{\zeta}_{2n}(z)$ to both sides of (\ref{e202br}), we get 
\begin{equation}
  \xi_{2n}(z)+\hat{\zeta}_{2n}(z)=2\hat{\zeta}_{2n}(z)-2^{1-z}\hat{\zeta}_n(z), 
\label{e203r}
\end{equation}
where the left-hand side is an order of $O(n^{-(1+\Re z)})$ for $n\to\infty$ 
shown in Appendix, whereas an order of each term is $O(n^{-\Re z})$.

  When we put $z=\rho$ which is one of the non-trivial zeros for the Riemann zeta function 
in (\ref{e203r}), 
take the absolute values and 
use the property of the zeta function that $1-\rho$ is also a zero as $\rho$ is, 
we get
\begin{eqnarray}
  2|\hat{\zeta}_{2n}(\rho)|&=&|2^{1-\rho}||\hat{\zeta}_n(\rho)|+O(n^{-(\sigma+1)})
\label{e2041r}\\
  2|\hat{\zeta}_{2n}(1-\rho)|&=&|2^{\rho}||\hat{\zeta}_n(1-\rho)|+O(n^{\sigma-2}),
\label{e205r}
\end{eqnarray}
where $\sigma=\Re \rho$.  

  Combining (\ref{e2041r}) with (\ref{e205r}), we get 
\begin{equation}
  |\hat{H}_{2n}(\rho)|=|2^{1-2\rho}||\hat{H}_n(\rho)|+O(n^{-2\sigma}). 
\label{e206r}
\end{equation}

  When we think about the limit of $n\to\infty$ in (\ref{e206r}), 
the left hand side will coincide with $|\hat{H}(\rho)|$ 
and $\displaystyle\lim_{n\to\infty}|\hat{H}_{n}(\rho)|$ converges to 
the same value in the right hand side. 
The reason why we have introduced the absolute values in (\ref{e2041r}) to (\ref{e206r}) is as follows. 
A dependence of $n$ in $\hat{H}_{n}(\rho)$ also appears as a function of $n$ in the argument. 
To eliminate this dependence, we have taken absolute values. 
After all we can conclude that the term $|2^{1-2\rho}|$ is equal to one 
which means that a real part of the zero $\Re\rho$ is identical to one half. 

\vskip 5mm


  For a while, we think about each value of $n$-dependent function for 
the non trivial zeros. 
We can easily get the $n$-dependent value 
and the next leading order for the $\hat{H}_n(z)$ as follows: 
\begin{eqnarray}
  \hat{\zeta}_n(\rho)&=&\frac{1}{2n^\rho}+O(n^{-(\sigma+1)}),\\
  \hat{\zeta}_n(1-\rho)&=&\frac{1}{2n^{1-\rho}}+O(n^{\sigma-2}).
\label{e200a}
\end{eqnarray}
So we can reach the relation
\begin{equation}
  \hat{H}_n(\rho)\equiv\frac{\hat{\zeta}_n(\rho)}{\hat{\zeta}_n(1-\rho)}
                 =n^{1-2\rho}+O(n^{-2\sigma}).
\label{e200-1a}
\end{equation}
When we think about the limit of $n\to\infty$ in (\ref{e200-1a}), 
the left hand side will converges the finite value which coincides with 
$\hat{H}(\rho)$ beside the argument. 
The real part of $1-2\rho$ must be equal to one in the right-hand side 
to have a finite value. 
Then we can conclude that a real part of the zero $\Re\rho$ is identical to one half.

\vskip 5mm
  For the use of non-regularized quantities $\zeta_n(z)$ appeared on right-hand side in Eq.(\ref{e123r}), 
the functional equation for $n$ can be also defined as follows:
\begin{equation}
  \zeta_n(z)\equiv H_n(z)\zeta_n(1-z),
\label{e201a}
\end{equation}
where $\zeta_n(z)$ is defined by Eq.(\ref{e101r}) and we refer $H_n(z)$ to the non-regularized coefficient as
\begin{eqnarray}
  H_n(z)&=&\frac{\zeta_n(z)}{\zeta_n(1-z)}=\frac{\hat{\zeta}_n(z)+n^{1-z}/(1-z)}{\hat{\zeta}_n(1-z)+n^z/z}\\
             &=&n^{1-2z}\left\{\frac{\hat{\zeta}_n(z)/n^{1-z}+1/(1-z)}{\hat{\zeta}_n(1-z)/n^z+1/z}\right\}\\
             &=&n^{1-2z}\left\{\frac{\hat{\zeta}_n(z)/n^{1-z}+1/(1-z)}{\hat{H}_n(z)^{-1}\hat{\zeta}_n(z)/n^z+1/z}\right\}
\label{e202a}
\end{eqnarray}
Considering the limit of $n\to\infty$ for $z=\rho$, 
we get a relation 
\begin{eqnarray}
  H_n(\rho)&=&n^{1-2\rho}\left\{\frac{\hat{\zeta}_n(\rho)/n^{1-\rho}+1/(1-\rho)}
                              {\hat{H}_n(\rho)^{-1}\hat{\zeta}_n(\rho)/n^\rho+1/\rho}\right\}\\
           &=&\frac{\rho}{1-\rho}n^{1-2\rho}+O(n^{-2\sigma}).
\label{e203}
\end{eqnarray}

  When we deal with the Euler's alternating series in Eq.(\ref{e102r}) for the Riemann zeta function, 
we can evaluate the function for $z$ even in the critical strip $0<\Re z<1$ as mentioned above. 
  The discussion about the ratio of $H_{2n}(\rho)$ and $H_n(\rho)$ is 
parallel to the regularized quantities\cite{Fujimoto5} which we have briefly reviewed above, 
by using Eq.(\ref{e202ar}), we get 
\begin{equation}
  \lim_{n\to\infty}{\zeta}_{2n}(\rho)=2^{1-\rho}\lim_{n\to\infty}{\zeta}_n(\rho)
\label{e207a}
\end{equation}
and using the property that $1-\rho$ is also a zero as well as $\rho$ is, we also get
\begin{equation}
  \lim_{n\to\infty}{\zeta}_{2n}(1-\rho)=2^{\rho}\lim_{n\to\infty}{\zeta}_n(1-\rho).
\label{e208a}
\end{equation}
  Combining Eqs.(\ref{e207a}) with (\ref{e208a}), we get 
\begin{equation}
  \lim_{n\to\infty}{H}_{2n}(\rho)=2^{1-2\rho}\lim_{n\to\infty}{H}_n(\rho). 
\label{e209a}
\end{equation}

  So the way to reach the conclusion is taking an absolute value 
in (\ref{e209a}) with (\ref{e206r})
\begin{equation}
  |2^{1-2\rho}|=\lim_{n\to\infty}\left|\frac{{H}_{2n}(\rho)}{{H}_n(\rho)}\right|
               =\lim_{n\to\infty}\left|\frac{\hat{H}_{2n}(\rho)}{\hat{H}_n(\rho)}\right|
               =\frac{|\hat{H}(\rho)|}{|\hat{H}(\rho)|}=1,
\label{e212a}
\end{equation}
where we have to use the regularized quantities because we do not know 
whether the non-regularized quantities limit to the same absolute value, 
but we can confirm that this is the fact and 
which claims that a real part of the zero $\Re\rho$ is equal to 
$\displaystyle\frac{1}{2}$. 


\vskip 5mm

  On the other hand, by using Eq.(\ref{e202ar}) and the property of the Euler's alternating series, 
we can get 
\begin{equation}
  \zeta_{2n}(\rho)=2^{1-\rho}\zeta_n(\rho)+O(n^{-\sigma}),
\label{e221a}
\end{equation}
and the $n$-dependent zeta function in each side can be reduced to $n$-th power using 
Eqs.(\ref{e123r}) and (\ref{e201a}), 
and we get the relation 
\begin{equation}
  H_{2n}(\rho)(2n)^{\rho}/\rho
                              =(2n)^{1-\rho}/(1-\rho)+O(n^{-\sigma})
\end{equation}
leads us to
\begin{equation}
  H_{2n}(\rho)=\frac{\rho}{1-\rho}(2n)^{1-2\rho}+O(n^{-2\sigma}),
\label{e222a}
\end{equation}
which is consistent with (\ref{e203}). 
In Eq.(\ref{e222a}) we can find the fact that even the non-regularized quantities $H_n(\rho)$ converges as
\begin{equation}
  \lim_{n\to\infty}\frac{H_n(\rho)}{n^{1-2\rho}}=\frac{\rho}{1-\rho}, 
\label{e223a}
\end{equation}
as far as a real part of the zero is equal to $\displaystyle\frac{1}{2}$, 
which is consistent with the Riemann hypothesis. 

  We can also show the relations concerning to regularized quantity for $n\to\infty$ as
\begin{equation}
  \lim_{n\to\infty}\frac{\hat{H}_n(\rho)}{n^{1-2\rho}}=1\ \ \ {\rm and}\ \ \ 
  \lim_{n\to\infty}\frac{\hat{\zeta}'_n(\rho)}{\hat{\zeta}'_n(1-\rho)}=-\hat{H}(\rho), 
\label{e224a}
\end{equation}
where the prime means the derivative for $z$.

  All these order estimation above are made use of the relation derived from Hardy and Littlewood\cite{Hardy}
\begin{equation}
  \hat{\zeta}(z)=\hat{\zeta}_n(z)+O(n^{-\Re z})
\end{equation}
for $|\Im z|\le 2\pi n/C$, where $C$ is constant greater than one, 
which means that $\Im z$ can be taken as large as $n$.

  Finally we have to mention that the regularized form for the Riemann zeta function above 
in the region $\frac{1}{2}\le\Re z<1$ coincides with the analytic continuation.

\vskip 5mm
\newpage
\renewcommand{\theequation}{\Alph{section}\arabic{equation}}
\setcounter{section}{1}
\setcounter{equation}{0}
\section*{Appendix}

\hspace{\parindent}
  Here we briefly show that the function of $z$ with a parameter $n$, 
namely the left-hand side of (\ref{e203r}) 
\begin{equation}
  h_{2n}(z)\equiv \xi_{2n}(z)+\hat{\zeta}_{2n}(z)
\label{a001}
\end{equation}
converges to zero rapidly compared with $\xi_{2n}(z)$ or $\hat{\zeta}_{2n}(z)$ 
for $n\to\infty$ at $z=\rho$, 
one of the non-trivial zeroes for the Riemann zeta function.
The definitions given in (\ref{e101r}) and (\ref{e123r}) are 
\begin{eqnarray}
  \zeta_n(z)&\equiv& \sum_{k=1}^n\frac{1}{k^z}\nonumber\\
  \hat{\zeta}_n(z)&\equiv& \zeta_n(z)-\frac{n^{1-z}}{1-z}\nonumber\\
  \hat{\zeta}(z)&\equiv& \lim_{n\to\infty}\hat{\zeta}_n(z).\nonumber
\end{eqnarray}

By using the analytical continuation of the Euler-Maclaurin sum formula, we can write 
\begin{equation}
  \hat{\zeta}(z)=\zeta_n(z)-\frac{n^{1-z}}{1-z}-\frac{1}{2n^z}+R_n(z),
\label{a002}
\end{equation}
where $R_n(z)$ represents a residue term including the Bernoulli terms, 
 which order of $o(n^{-\Re z})$ for $n\to\infty$.
We write down
\begin{equation}
  \hat{\zeta}(z)-\hat{\zeta}_n(z)=-\frac{1}{2n^z}+R_n(z).
\label{a003}
\end{equation}

On the other hand, the definition for the Euler alternating series 
given in (\ref{e102r}) is 
\begin{eqnarray}
  \xi_n(z)&\equiv& \sum_{k=1}^n\frac{(-1)^{k-1}}{k^z}\nonumber\\
  \xi(z)&\equiv& \lim_{n\to\infty}\xi_n(z).\nonumber
\end{eqnarray}
This immediately means that
\begin{equation}
  \xi(z)=(1-2^{1-z})\hat{\zeta}(z),
\label{a004}
\end{equation}
\begin{equation}
  \xi_{2n}(z)
             =\hat{\zeta}_{2n}(z)-2^{1-z}\hat{\zeta}_n(z)
\label{a005}
\end{equation}
and
\begin{equation}
  h_{2n}(z)
           = 2\hat{\zeta}_{2n}(z)-2^{1-z}\hat{\zeta}_n(z).
\label{a006}
\end{equation}

Here we note the difference of first order as
\begin{eqnarray}
  g_{2n}(z)&\equiv& \frac{\xi_{2n-1}(z)+\xi_{2n}(z)}{2}\nonumber\\
           &=& \xi_{2n}(z)+\frac{1}{2(2n)^z}
\label{a007}
\end{eqnarray}
and we get 
\begin{eqnarray}
  & & \hat{\zeta}(z)-h_{2n}(z)\nonumber\\
  &=& \hat{\zeta}(z)-\{2\hat{\zeta}_{2n}(z)-2^{1-z}\hat{\zeta}_n(z)\}\nonumber\\
  &=& \{\hat{\zeta}(z)-\hat{\zeta}_{2n}(z)\}-\{\hat{\zeta}_{2n}(z)-2^{1-z}\hat{\zeta}_n(z)\}\nonumber\\
  &=& \left\{-\frac{1}{2(2n)^z}+R_{2n}(z)\right\}-\xi_{2n}(z)\nonumber\\
  &=& -g_{2n}(z)+R_{2n}(z).
\label{a008}
\end{eqnarray}
Putting $z=\rho$, one of non-trivial zeroes, we get a relation 
\begin{equation}
  h_{2n}(\rho)=g_{2n}(\rho)-R_{2n}(\rho).
\label{a009}
\end{equation}

Meanwhile we evaluate $g_{2n}(z)$ as
\begin{eqnarray}
  g_{2n}(z)&=& \xi_{2n}(z)+\frac{1}{2(2n)^z}\nonumber\\
           &=& \zeta_{2n}(z)-2^{1-z}\zeta_n(z)+\frac{1}{2(2n)^z}\nonumber\\
           &=& \left\{\hat{\zeta}(z)+\frac{(2n)^{1-z}}{1-z}+\frac{1}{2(2n)^z}-R_{2n}(z)\right\}\nonumber\\
            && -2^{1-z}\left\{\hat{\zeta}(z)+\frac{n^{1-z}}{1-z}+\frac{1}{2n^z}-R_n(z)\right\}
               +\frac{1}{2(2n)^z}\nonumber\\
           &=& (1-2^{1-z})\hat{\zeta}(z)-R_{2n}(z)+2^{1-z}R_n(z)\nonumber\\
           &=& \xi(z)-R_{2n}(z)+2^{1-z}R_n(z)
\label{a010}
\end{eqnarray}
and again putting $z=\rho$, we get another relation
\begin{equation}
  g_{2n}(\rho)=-R_{2n}(\rho)+2^{1-\rho}R_n(\rho).
\label{a011}
\end{equation}

By using (\ref{a009}) and (\ref{a011}), we get 
\begin{equation}
  h_{2n}(\rho)=-2R_{2n}(\rho)+2^{1-\rho}R_n(\rho).
\label{a013}
\end{equation}

Eqs.(\ref{a011}) and (\ref{a013}) are both order of $o(n^{-1/2})$, which shows 
that the $h_{2n}(\rho)$ converges to zero more rapid than $\xi_{2n}(\rho)$ or $\hat{\zeta}_{2n}(\rho)$.

%

\vskip 5mm


\noindent
\begin{thebibliography}{(00)}

\bibitem{Conrey} J.B. Conrey,  {\it The Riemann hypothesis}, Notice of the AMS 50 (2003), 341-353.
\bibitem{Fujimoto1} M. Fujimoto and K. Uehara, {\it Regularization for zeta functions 
with physical applications I}, arXiv:math-ph/0609013 (2006).
\bibitem{Fujimoto2} M. Fujimoto and K. Uehara, {\it Regularization for zeta functions 
with physical applications II}, arXiv:math-ph/0702011 (2007).
\bibitem{Fujimoto3} M. Fujimoto and K. Uehara, {\it Regularizations of the Euler product 
representation for zeta functions and the Birch--Swinnerton-Dyer conjecture}, 
arXiv:math-ph/0709.0762 (2007).
\bibitem{Fujimoto4} M. Fujimoto and K. Uehara, {\it Regularized Euler product for 
the zeta function and the Birch and Swinnerton-Dyer and the Beilinson conjecture}, 
arXiv:math-ph/0811.2644 (2008).
\bibitem{Fujimoto5} M. Fujimoto and K. Uehara, {\it A Brief Note on the Riemann hypothesis}, 
arXiv:math-ph/0906.1099 (2009), math-ph/0906.1099v2 (2011).
\bibitem{Guy} R. K. Guy, Products Taken over Primes \S B87 in 
{\it Unsolved Problems in Number Theory}, 2nd ed. New York: Springer-Verlag, pp. 102-103, 1994. 
\bibitem{Koblitz} N. Koblitz, {\it p-adic numbers, p-adic analysis and zeta function}, 
Graduate Texts in Mathematics, 58, 2nd edition, Springer, 1998.
\bibitem{Kunimasa} T. Kunimasa and K. Uehara, Nucl. Phys. B279 (1987)~608-640.
\bibitem{Leeuwen} J.M.J. van Leeuwen, 
  Proc. Kon. Ned. Acad. v. Wetenschap Centenary issue 100 (1997)~57-63.
\bibitem{Sarnak} Z. Rudnick and P. Sarnak, C.R. Acad. Sci. Paris 319 (1994)~1027-1032.
\bibitem{Hardy} G.H. Hardy and J.E. Littlewood, Proc. London Math. Soc., 
{\bf s2-29}(1) (1929),~81-97.

\end{thebibliography}
\end{document}